\documentclass[a4paper,12pt]{amsart}
\usepackage{amssymb}
\usepackage[]{latexsym}

\usepackage{amsthm}

\usepackage{lscape}

\usepackage[latin1]{inputenc}
\usepackage{color}
\usepackage[all,dvips]{xy}

\usepackage[dvips]{graphicx}

\newtheorem{thm}{Theorem}
\newtheorem{cor}[thm]{Corollary}
\newtheorem{lem}[thm]{Lemma}

\newtheorem{prop}[thm]{Proposition}

\theoremstyle{remark}
\newtheorem{rem}[thm]{Remark}

\theoremstyle{definition}

\newenvironment{pf}{\par\noindent{\bf Proof.}\enspace\ignorespaces}{\qed\par\par}

\newenvironment{pfm}{\par\noindent{\bf Proof of Theorem \ref{Mthm}.}\enspace\ignorespaces}{\qed\par\par}

\def\qed{\hfill $\Box$}

\newcommand{\pmodd}[1]{\,({\rm mod\,}{#1}) }

\newcommand{\bQ}{{\mathbb{Q}}}
\newcommand{\bF}{{\mathbb{F}}}

\newcommand{\bP}{{\mathbb{P}}}
\newcommand{\bZ}{{\mathbb{Z}}}
\newcommand{\bA}{{\mathbb{A}}}

\newcommand{\cO}{{\mathcal{O}}}
\newcommand{\cJ}{{\mathcal{J}}}

\DeclareMathOperator{\Jac}{Jac}

\title{On symmetric square values of quadratic polynomials}

\author[Enrique Gonz\'alez-Jim\'enez]{Enrique Gonz\'alez-Jim\'enez}
\address{Universidad Aut{\'o}noma de Madrid, Departamento de Matem{\'a}ticas and Instituto de Ciencias Matem{\'a}ticas (CSIC-UAM-UC3M-UCM), Madrid, Spain}
\email{enrique.gonzalez.jimenez@uam.es}
\thanks{The first author was supported in part by grants MTM 2009-07291 (Ministerio de Educaci{\'o}n y Ciencia, Spain)
and  CCG08-UAM/ESP-3906 (Universidad Auton{\'o}ma de Madrid-Comunidad
de Madrid, Spain). The second author was partially supported by
the grant MTM2009-10359.}

\author[Xavier Xarles]{Xavier Xarles}
\address{Departament de Matem\`atiques\\Universitat Aut\`onoma de
Barcelona\\08193 Bellaterra, Barcelona, Spain}
\email{xarles@mat.uab.cat}


\date{\today}

\begin{document}
\maketitle

\begin{abstract}
We prove that there does not exist a non-square quadratic
polynomial with integer coefficients and an axis of symmetry which
takes square values for $N$ consecutive integers for $N=7$ or
$N\geq 9$. At the opposite, if $N\le 6$ or $N=8$ there are
infinitely many.
\end{abstract}


\section{Introduction}

In this note we are dealing with the following problem. Given a
degree two polynomial $f(x)=ax^2+bx+c \in \bZ[x]$ which is not a
square of a degree one polynomial, how many consecutive integers
values $f(i)$ can be squares in $\bZ$? This problem has been
considered by D. Allison in \cite{All1} and \cite{All2}, who found
infinitely many examples with eight consecutive values, and by A.
Bremner in \cite{Bre}, who found more examples with seven
consecutive values.

The examples found by Allison are all  by polynomials
which are symmetric with an axis of symmetry midway between two
integers. This means that, after some easy translation, all the
examples are of the form $f(x)=a(x^2+x)+c$ and the values are
$f(i)$ for $i=-3,-2,-1,0,1,2,3$ and $4$. This result was obtained
by translating the problem to computing rational points on some
elliptic curve which has rank one.

On the other hand, Bremner in \cite{Bre} shows that there does not
exists any example which is symmetric with an axis of symmetry about
an integral value and with 7 values, by showing that these examples
would be described by rational points in some rank zero elliptic
curve, which has 12 points, all corresponding to the polynomial
$f(x)$ being the square of a polynomial.

In the same paper, Bremner asks if there are examples as the ones
found by Allison, but with 10 consecutive squares. The problem
translates to finding all the rational points of a genus 5 curve,
a fact already noticed by Allison and by Bremner. He conjectures
that there should be no such example.

In this note we prove this conjecture, and so, combining with the
results of Bremner and Allison, we get the following theorem.

\begin{thm}\label{Mthm}
Let $N$ be a positive integer and $\mathcal{B}_N$ the set
consisting of non-square quadratic polynomials $f(x)=a x^2+b x+c
\in \bZ[x]$ that takes square values for $N$ consecutive integers
values of x, $r$,$r+1$,...,$r+N-1$ and $f(r)=f(r+N-1)$, for some
$r\in\bZ$. Then
$$
\#\mathcal{B}_N=\left\{
\begin{array}{ccl}
\infty & & \mbox{if $N\le 6$ or $N=8$},\\[1mm]
0 & & \mbox{if $N=7$ or $N\ge 9$}.
\end{array}
\right.
$$
\end{thm}



To show this result we will use similar techniques as the one we
use in \cite{GX} to study the arithmetic progressions of squares
over quadratic fields. In fact, the problem we study here is in
some sense a generalization to higher dimension of the old result
by Fermat about arithmetic progressions of squares, and the
problems in \cite{GX} and in \cite{X} are generalizations to
higher degree (in the sense of the field of numbers involved).
However, we can say almost nothing about if there exists a maximum
number of consecutive square values taken by a non-square
quadratic polinomial, or even if this number is $8$, as the known
examples suggest.

\section{Translation to Geometry}\label{sec2}

Fix a polynomial $f(x)=a x^2+b x+c \in \bZ[x]$ which takes square
values for $N$ consecutive integers values of $x$, of the form
$r$,$r+1$,...,$r+N-1$, and $f(r)=f(r+N-1)$. Suppose that $N$ is even.
Then, after translation by $-r-N/2$, we can suppose that $r=-N/2$.
We get then that $f(x)$ has the form $f(x)=a(x^2+x)+c$, and we are
asking to have $f(i)=x_i^2$ for $i=0,\dots,N/2-1$ and $x_i\in \bZ$.

Now, suppose that $N=10$. The conditions we get from $f(i)=x_i^2$
for $i=0,\dots,4$ are given by the following equations:
$$
C:\left\{
\begin{array}{rcrcc}
2 x_0^2 \!\!\! & - &\!\!\!\!\! 3 x_1^2 &\!\!\! + &\!\!\! x_2^2=0,\\[1mm]
5x_0^2 \!\!\!  & - &\!\!\!\!\! 6 x_1^2 &\!\!\! + &\!\!\! x_3^2=0,\\[1mm]
9 x_0^2 \!\!\! & - &\!\!\!\!\! 10x_1^2 &\!\!\! + &\!\!\! x_4^2=0,
\end{array}
\right.$$
which determine a genus 5 curve $C$ in $\bP^4$. Any point
$P:=[x_0:x_1:x_2:x_3:x_4]$ of this curve defined over $\bQ$ will
give us a polynomial $f(x)$ as before, by setting $c=x_0^2$ and
$a=(x_1^2-x_0^2)/2$. Observe that the pairs $(a,c)$ are well
defined modulo multiplication by a square number, which will
produce the same polynomial but multiplied by a square number, a
case that we consider equivalent.

Now, the solutions given by $P=[\pm 1:\pm 1:\pm 1:\pm 1:\pm 1]$
correspond to the case $a=0$, so the polynomial is in fact
constant. There are also the solutions given by $P=[\pm1:\pm 3:\pm
5:\pm 7:\pm 9]$, which correspond to the case $a=4$ and $c=1$, so
$f(x)=(2x+1)^2$. Our aim will be to show that these are the only
rational points.

First of all, observe that the curve $C$ has degree $2$ maps
$\Phi_n$ to five distinct elliptic curves $E_n$, for $n=0,1,2,3,4$.
They can be described easily as intersection of two quadrics in
$\bP^3$, by taking first two of the three quadrics describing $C$,
which give three of them, and transforming the equations in order
to get more quadratic forms involving only three variables, which
gives the other two. We denote by $F_n$ to the genus one curve obtained by removing the variable $x_n$ for $n=0,1,2,3,4$.
We have obtained the following equations for such a curves:
$$
\begin{array}{lrlll}
F_0\,:&x_{1}^2 & =16\,t^4 -144\,t^3+340\,t^2-252\,t+49\,, & \quad & t=\frac{x_4+x_2}{x_4-x_3},\\[1.5mm]
F_1\,:&x_{0}^2 & =16\,t^4 -160\,t^3+384\,t^2-280\,t+49\,, & \quad & t=\frac{x_4+x_2}{x_4-x_3},\\[1.5mm]
F_2\,:& x_4^2 & =36\,t^4 + 96\,t^3 - 236\,t^2 + 80\,t + 25\,, & \quad & t=\frac{x_3+x_1}{x_3-x_0},\\[1.5mm]
F_3\,:&x_{1}^2 & =100\,t^4-360\,t^3+472\,t^2-252\,t+49\,, & \quad & t=\frac{x_4+x_2}{x_4-x_0},\\[1.5mm]
F_4\,:&x_2^2 & =36\,t^4 - 72\,t^3 + 72\,t^2 - 60\,t + 25\,, & \quad &t= \frac{x_3+x_1}{x_3-x_0}.\\[1.5mm]
\end{array}
$$
The second column is in fact equivalent to the forgetful map
$\varrho_n:C\longrightarrow F_n$. Observe that there is always a
quadratic form involving only three of the variables $x_n$, for
any choice of them. The curve $F_n$ is isomorphic to the one given
by equations not involving the variable $x_n$.

All these genus one curves have rational points over $\bQ$,
therefore they are elliptic curves over $\bQ$. A Weierstrass model
of the quartic $F_n$ is denoted by $E_n$, $n=0,1,2,3,4$:
$$
\begin{array}{l}
E_0\,:y^2 = x(x-8)(x+27),\\[1.1mm]
E_1\,:y^2 =  x(x-12)(x+30), \\[1.1mm]
E_2\,:y^2 =  x(x+4)(x+54), \\[1.1mm]
E_3\,:y^2 =   x(x-7)(x+20), \\[1.1mm]
E_4\,:y^2 =   x(x-12)(x-15). \\[1.1mm]
\end{array}
$$
Using the labeling of the Cremona's tables \cite{cremonaweb}, one
can check that $E_0= $ \verb+1680G2+, $E_1= $ \verb+20160BG2+,
$E_2=$ \verb+960H2+$, E_3= $ \verb+840H2+ and $E_4=$ \verb+360E2+.

So, if one of such elliptic curves have a finite number of
rational points, then the problem of computing $C(\bQ)$ becomes
easy. Now, it is a straightforward computation to check that the
torsion subgroup of $E_n(\bQ)$ is isomorphic to
$\bZ/2\bZ\oplus\bZ/2\bZ$. On the other hand, the set $\Phi_n([\pm
1:\pm 1:\pm 1:\pm 1:\pm 1])$ has cardinality eight. Therefore the
rank of  $E_n(\bQ)$ is greater than one, that is, we cannot use
these argument to determine $C(\bQ)$.  In fact, we can easily
compute by descent (or, better, using some algebraic computational
system like \verb+Magma+ or \verb+Sage+, or even better, checking
at the Cremona's tables), that $\mbox{rank}_\bZ\, E_n(\bQ)=1$ for
$n\ne 1$ and $\mbox{rank}_\bZ\, E_1(\bQ)=2$.

Another approach is the Chabauty's method. This could be used when
the Jacobian of the curve has rank less than the genus of the
curve. But in this case the rank  of $\Jac(C)$ is greater than
$5$, the genus of $C$, since  $\Jac(C)$ is $\bQ$-isogenous to the
product of $E_n$, for $n=0,1,2,3,4$. So we cannot apply this
method. Other methods, like the Manin-Drinfeld's method, cannot be
applied either.

\section{Two descent and Covering Collections}

In order to actually compute the rational points in the curve $C$,
we will apply the covering collections technique, as developed by
Coombes and Grant \cite{CG}, Wetherell \cite{We} and others, and
specifically a modification of what is now called the elliptic
Chabauty method developed by  Flynn and Wetherell in \cite{FW} and
by Bruin in \cite{Br}.

The method has two parts. Suppose we have a curve $C$ over a
number field $K$ and an unramified map $\chi:C'\to C$ of degree
greater than one and may be defined over a finite extension $L$ of
$K$. We consider the distinct unramified coverings
$\chi^{(s)}:C'^{(s)}\to C$ formed by twists of the given one, and
we get that
$$C(K)=\bigcup_{s} \chi^{(s)}(\{P\in C'^{(s)}(L)\ : \ \chi^{(s)}(P)\in C(K)\}) ,$$
the union being disjoint. Only a finite number of twists do have
rational points, and the finite (larger) set of twists having
points locally everywhere can be explicitly described. The first
part is to compute this set of twists, and the second to compute
the points $ P\in C'^{(s)}(L)$ such that $\chi^{(s)}(P)\in C(K)$.
This second part depends on having a nice quotient of the curves
$C'^{(s)}$, for example a genus one quotient, where it is possible
to do the computations. In this section we will concentrate on the
first part.

The coverings we are going to consider are Galois coverings with
Galois group isomorphic to $(\bZ/2\bZ)^2$. This coverings are in
principle easy to construct. One only needs to have an isogeny map
from an abelian variety $A$ to the jacobian $\Jac(C)$ of the curve
$C$ with kernel isomorphic, as group scheme, to the group
$(\bZ/2\bZ)^2$. Since in our case the jacobian $\Jac(C)$ is
isogenous to a product of elliptic curves $E_i$, they can be
constructed by choosing two such elliptic curves and one degree
two isogeny in each of them.

Moreover, the elliptic curves $E_i$ have all they $2$-torsion
points defined over $\bQ$, hence the coverings we are searching
for will be defined over $\bQ$. On the other hand, the genus one
quotients of such coverings that we will use in the next section
are, in general, not defined over $\bQ$, but in a quadratic or in
a biquadratic extension. The way we will construct the coverings,
by using a factorization of the quartic polynomials, will give us
directly also the genus one quotient and the field where it is
defined.

In order now to construct the coverings of the curve $C$, we first
rewrite the equations of the curve in the following form:
$$
C\,:\, \left\{
\begin{array}{rcl}
y^2&\!\!\!=&\!\!\!  q(t)=36\,t^4 - 72\,t^3 + 72\,t^2 - 60\,t + 25,   \\[1mm]
z^2 &\!\!\! =&\!\!\!   p(t)=36\,t^4 + 96\,t^3 - 236\,t^2 + 80\,t +25.
\end{array}\right.
$$
This model has two natural maps to the genus one curves
$F_4:y^2=q(t)$ and $F_2:z^2=p(t)$, whose jacobians are $E_4$ and
$E_2$ respectively.

First, we concentrate on the unramified degree two coverings of
the genus one curve given by a quartic model. If $F$ is such a
genus one curve defined over a field $K$, given by the equation
$y^2=r_1(x)r_2(x)$, where $r_1(x)$ and $r_2(x)$ are degree two
polynomials defined over an extension $L$ of $K$, we consider the
degree two unramified covering $\chi:F'\to F$ with affine part in
$\bA^3$ given by the zeros of the polynomials $y_1^2=r_1(x)$ and
$y_2^2=r_2(x)$, the map given by $\chi(x,y_1,y_2)=(x,y_1y_2)$.
For any $\delta \in L^*$, we consider the curve $F'^{(\delta)}$
given by the equations $\delta y_1^2=p_1(x)$ and $\delta
y_2^2=p_2(x)$, and the map to $F$ defined by
$\chi^{(\delta)}(x,y_1,y_2)=(x,y_1y_2/\delta^2)$. Then
$F'^{(\delta)}$ are all the quadratic twists of $F'$, and there
exists a finite set $\Delta_L(\chi)\subset L^*$ such that
$$F(K) \subseteq \bigcup_{\delta \in \Delta_L(\chi)} \chi^{(\delta)}(\{(x,y_1,y_2)\in F'^{(\delta)}(L)\ : \ x\in K  \mbox{ or } x=\infty \}).$$

First, we consider the case of $F_2$, given as $y^2=p(t)$, where
$$p(t)=36\,t^4 + 96\,t^3 - 236\,t^2 + 80\,t +
25=(6\,t^2-4\,t-1\,)(6\,t^2+20\,t-25\,).$$

\begin{lem}\label{deltaF2}
Consider the degree two covering defined over $\bQ$ given by
$$
F_2'\,:\, \left\{
\begin{array}{rcl}
z_1^2&\!\!\!=&\!\!\!  p_1(t)=6\,t^2-4\,t-1,   \\[1mm]
z_2^2 &\!\!\! =&\!\!\!   p_2(t)= 6\,t^2+20\,t-25\,,
\end{array}\right.
$$
together with the natural map $\psi_2:F_2'\to F_2$ given by
$\psi_2(t,z_1,z_2)=(t,z_1z_2)$. Then $\Delta_{\bQ}(\psi_2):=\{ \pm
1,\pm 6 \}$, hence
$$F_2(\bQ) \subseteq \bigcup_{\delta \in \{\pm 1, \pm 6 \}}
\psi_2^{(\delta)}(\{(t,z_1,z_2)\in F_2'^{(\delta)}(\bQ)\}).$$
\end{lem}

\begin{pf} It is easy to show and well-known that $\psi_2$ is an
unramified degree two covering of $F_2$. Since $F_2'(\bQ)\ne
\emptyset$, because it contains the point $P':=(1,1,1)$, we can
identify $F_2'$ with an elliptic curve $E'_2$, by sending this
point $P'$ to $O$, and, identify $F_2$ with the elliptic curve
$E_2$ by sending the point $P:=\chi(P')=(1,1)$ to the point $O$.
We get then an unramified degree two covering $\phi_2:E_2'\to
E_2$, which must be a degree $2$ isogeny. With appropriate choices
of the identifications, we can get this isogeny in the standard
form (see, for example, \cite[III.4.5]{Sil} or \cite[\S 8.2]{Coh}
). After some computations we get the map $\phi_2 : E'_2 \to E_2$
defined by
$$
\phi_2(x,y)=\left(\frac{y^2}{4x^2},\frac{y(x^2-2500)}{8x^2}\right),
$$
where $E_2\,:y^2 =  x(x-4)(x+54)$ and $E_2'\,: y^2=x(x^2 - 116x +
2500)$.

Now, the quadratic twists $F_2'^{(\delta)}$ which locally have
rational points correspond to the elements of the Selmer group
$\mbox{Sel}(\phi_2)$. After identifying $\mbox{Sel}(\phi_2)$ with
a subgroup of $\bQ^*/(\bQ^*)^2$ in the standard way, the
identification sends $\delta\in \bQ^*$ to its class modulo
squares. A standard $2$-descent calculation gives that
$\mbox{Sel}(\phi_2)=\{ \pm 1, \pm 6 \}$. But now, by using that
 $E_2(\bQ)$ contains the points $(0,0)$, $(-54, 0)$ and $(36,
 -360)$, one can see that
 all the elements of the Selmer group $\mbox{Sel}(\phi_2)$ correspond
 to elements of $E(\bQ)$. These are exactly the $\delta$'s
 such that $F_2'^{(\delta)}(\bQ)\ne \emptyset$.
\end{pf}

\

Now we consider the curve $F_4$, given as $y^2=q(t)$, where
$q(t)=36\,t^4 - 72\,t^3 + 72\,t^2 - 60\,t + 25 $. Observe that the
polynomial $q(t)$ is irreducible over $\bQ$, but it factorizes
over some quadratic extensions as product of two degree 2
polynomials. Over $\bQ(\sqrt{6})$ we have
$$
q(t)=(5+2\sqrt{6})(6\,t^2+2\sqrt{6}\,t+5-6\sqrt{6})\cdot (5-2\sqrt{6})\,(6\,t^2-2\sqrt{6}\,t+5+6\sqrt{6}).
$$
\begin{lem}\label{deltaF4}
Consider the degree two covering defined over $\bQ$ given by
$$
F_4'\,:\, \left\{
\begin{array}{rcl}
y_1^2&\!\!\!=&\!\!\!  q_1(t)=(5+2\sqrt{6})\,(6\,t^2+2\sqrt{6}\,t+5-6\sqrt{6})\,  , \\[1mm]
y_2^2 &\!\!\! =&\!\!\!q_2(t)=
(5-2\sqrt{6})\,(6\,t^2-2\sqrt{6}\,t+5+6\sqrt{6}),
\end{array}\right.
$$
together with the natural map $\psi_4:F_4'\to F_4$ given by
$\psi_4(t,y_1,y_2)=(t,y_1y_2)$. Then
$\Delta_{\bQ(\sqrt{6})}(\psi_4):=\{1,2,5,10\}$, hence
$$F_4(\bQ) \subseteq \bigcup_{\delta \in \{1,2,5,10\}} \psi_4^{(\delta)}(\{(t,y_1,y_2)\in F_4'^{(\delta)}(L) \ :
\ t\in \bQ \mbox{ or } t=\infty \}).$$
\end{lem}

\begin{pf} As in the proof of the lemma above, observe that
$F_4'(\bQ(\sqrt{6}))$ contains the point $(1,1,1)$, such that
$\psi_4(1,1,1)=(1,1)\in F_4(\bQ)$. Then the degree two covering
$\psi_4:F_4'\to F_4$ defined over $\bQ(\sqrt{6})$ corresponds, by taking some
isomorphisms to the respective jacobians, to the $2$-isogeny
$\phi_4 : E'_4 \to E_4$ defined by
$$
\phi_4(x,y)=\left(\frac{y^2}{4x^2},\frac{y(x^2-9)}{8x^2}\right),
$$
where $E_4: y^2=x(x-12)(x-15)$ and $E'_4: y^2=x(x^2 + 54x + 9)$,
which is the dual isogeny of the $2$-isogeny corresponding to the
$2$-torsion point $P=(0,0)\in E_4(\bQ)$. Now, a descent
computation shows that $\mbox{Sel}(\phi_4)=\{1, 3, 2, 6, 5, 15,
10, 30\}$. But observe now that two $\delta$ and $\delta'\in \bQ$
that are equivalent modulo squares over $\bQ(\sqrt{6})^*$ give isomorphic
coverings $\psi_4^{(\delta)}$. Hence we need only to consider the
set $\Delta_{\bQ(\sqrt{6})}(\psi_4)$ which is $\mbox{Sel}(\phi_4)$ modulo
$(\bQ(\sqrt{6})^*)^2$, which gives the result.
\end{pf}

\

We take now the unramified covering $\xi: C'\to C$ defined by the
equations:
$$
C'\,:\,
\left\{
\begin{array}{rcl}
y_1^2&\!\!\!= &\!\!\! q_1(t)=  (5+2\sqrt{6})\,(6\,t^2+2\sqrt{6}\,t+5-6\sqrt{6})\,,\\[1mm]
y_2^2&\!\!\!= &\!\!\! q_2(t)= (5+2\sqrt{6})\,(6\,t^2-2\sqrt{6}\,t+5+6\sqrt{6})\,, \\[1mm]
z_1^2&\!\!\!= &\!\!\! p_1(t)= 6\,t^2-4\,t-1\,, \\[1mm]
z_1^2&\!\!\!= &\!\!\! p_2(t)=  6\,t^2+20\,t-25\,,\\[1mm]
\end{array}\right.
$$
which is a curve of genus $17$.

The lemmata above computes the relevant twists to be consider.

\begin{cor}\label{deltaC} The set of revelant twists is equal to
$$\Delta:=\{ (\delta_2,\delta_4)\in \bQ(\sqrt{6})^* \ | \ \delta_i\in
\Delta_{\bQ(\sqrt{6})}(\phi_i),\  i=2, 4\},$$ where $\Delta_{\bQ(\sqrt{6})}(\phi_2)=\{ \pm 1\}$ and
$\Delta_{\bQ(\sqrt{6})}(\phi_4)=\{1,2,5,10\}$, which correspond to a set of
representatives in ${\bQ(\sqrt{6})}$ of the image of Selmer groups of $\phi_i$ ($i=2,4$)
in ${\bQ(\sqrt{6})}^*/({\bQ(\sqrt{6})}^*)^2$ via the natural maps.

Hence, $$C(\bQ) \subseteq \bigcup_{\delta \in \Delta}
\chi^{(\delta)}(\{(t,y_1,y_2,z_1,z_2)\in C'^{(\delta)}(\bQ(\sqrt{6}))\ : \ t\in \bQ
 \mbox{ or } t=\infty \}),$$ where $C'^{(\delta_2,\delta_4)}$ is the curve defined by
$$C'^{(\delta_2,\delta_4)}\,:\,\{\delta_4 y_1^2=q_1(t) \,\,,\,\,
\delta_4 y_2^2=q_2(t) \,\,,\,\, \delta_2 z_1^2=p_1(t) \,\,,\,\, \delta_2 z_2^2=p_2(t)\}. $$
\end{cor}

\begin{pf} After the lemmata above, we only need to observe that
$\Delta_{\bQ}(\phi_2)=\{ \pm 1,\pm 6\}$ becomes, after taking the
image in ${\bQ(\sqrt{6})}^*/({\bQ(\sqrt{6})}^*)^2$, the set
$\Delta_{\bQ(\sqrt{6})}(\phi_2)=\{ \pm 1\}$ .
\end{pf}

\

One can reduce even further the set of revelant twists to be
considered by using the natural automorphisms of $C$ given by
interchanging the sign of one of the coordinates (in the first
model of $C$).

\begin{cor}\label{deltaCinv} Let $\tau_i$ be the automorphisms of $C$ given by $\tau_i(x_i)=-x_i$,
and $\tau_i(x_j)=x_j$ if $j\ne i$, for $i=0,1,2,3$ and $4$, and
let $\Upsilon$ be the subgroup they generate. Let $\Delta'=\{ (1,1),(-1,1)\}$. Then, for any $P\in C(\bQ)$,
there exits $\tau \in \Upsilon$ and $\delta\in \Delta'$ such that
$$\tau(P) \in \chi^{(\delta)}(\{(t,y_1,y_2,z_1,z_2)\in C'^{(\delta)}(\bQ(\sqrt{6}))\ : \
t\in \bQ  \mbox{ or } t=\infty \}),$$ where $\chi^{(\pm 1
,1)}(x,y_1,y_2,z_1,z_2)=(x,y_1y_2,z_1z_2)$.
\end{cor}

\begin{pf}
In order to prove the result, it is enough to show that for any point
$P\in C(\bQ)$, and for any $\delta_4\in
\Delta_{\bQ(\sqrt{6})}(\phi_4)$, there exists $\delta_2\in
\Delta_{\bQ(\sqrt{6})}(\phi_2)$ and $\tau\in \Upsilon$ such that
$\tau(P)\in
\chi^{(\delta_2,\delta_4)}C'^{(\delta_2,\delta_4)}(\bQ(\sqrt{6}))$.
Therefore, one can reduce to show that for any  $\delta_4\in
\Delta_{\bQ(\sqrt{6})}(\phi_4)$, there exists $\tau\in \Upsilon$
such that the image of $\varrho_2(\tau(P))$ in
$\mbox{Sel}(\phi_2)$ is equal to $\delta_2$ modulo
$(\bQ(\sqrt{6})^*)^2$, where $\varrho_2: C\to F_2$ is the map
given in section \ref{sec2}.

Observe that the involutions $\tau_i$ for $i=0,1,3,4$ determine
involutions in $F_4$, which in turn determine involutions
$\tau_i'$ in $E_4$ once fixed the isomorphism between $F_4$ and
$E_4$. These involutions must be of the form $\tau_i'(Q)=-Q+Q_i$
for some $Q_i\in E(\bQ)$, since they do have fixed points. Hence,
the involutions $\tau_i'$ are determined once we know the image of
a single point $Q$. Thus, if we know the result for just one point
$P\in C(\bQ)$, we will obtain the result for all points in
$C(\bQ)$.

Take $P:=[1:1:1:1:1]$. Then one shows easily that the image of
$\varrho_2(\tau_i(P))$ in the Selmer group $\mbox{Sel}(\phi_2)$
for $i=0,1,3$, together with $\varrho_2(P)$, covers all the group,
which proves the result.
\end{pf}

\
\begin{rem} An easy computation by using the maps $\varrho_n: C\to
F_n$ given in section \ref{sec2}, shows that the involutions $\tau_i$ take
the following form in the model of $C$ given by $y^2=q(t)$ and
$z^2=p(t)$:
$$
\begin{array}{l}
\tau_0(t,y,z)=\left(\frac{6t-5}{6(t-1)},\frac{y}{6(t-1)^2},\frac{z}{6(t-1)^2}\right),\\[2mm]
\tau_1(t,y,z)=\left(\frac{5(t-1)}{6t-5},\frac{5y}{(6t-5)^2},\frac{5z}{(6t-5)^2}\right),\\[2mm]
\tau_3(t,y,z)=\left(\frac{5}{6t},\frac{5 y}{6t^2},\frac{5
z}{6t^2}\right),
\end{array}
$$
$\tau_2(t,y,z)=(t,y,-z)$ and $\tau_4(t,y,z)=(t,-y,z)$. These can
be used also to show the last Corollary.
\end{rem}

Observe that the known points in $C(\bQ)$, corresponding to the
points $[1:\pm1:\pm1:\pm1:\pm1]$ and $[1:\pm3:\pm5:\pm7:\pm9]$,
give rise to the points in $C'^{( 1,1)}(\bQ(\sqrt{6}))$ with $t=1$
and in $C'^{(-1,1)}(\bQ(\sqrt{6}))$ with $t=1/2$, respectively.

Now, in order to compute the points $(t,y_1,y_2,z_1,z_2)$ in
$C'^{(\pm 1,1)}(\bQ(\sqrt{6}))$ such that $t\in \bQ$, we consider
the following natural genus one quotients of $C'^{(\pm 1,1)}$
defined by
$$H_{i,j}^{(\pm)} \ : \ \pm w^2=q_i(t)p_j(t)$$
for $i,j=1,2$. We have four of them for any sign, corresponding in
fact to the factors of a natural genus $4$ quotient of any of the
curves $C'^{(\pm 1,1)}$, which is defined over $\bQ$.

Hence, we only need to compute, the points
$$\{(t,w) \in H_{(i,j)}^{(\pm 1,1)} \ : \ t\in \bQ \mbox{ or } t=\infty\}$$
for some $(i,j)$, and we are done. But this can be done by using
the elliptic Chabauty method.

 The following diagram illustrates all the curves
and morphisms involved in our problem:
$$
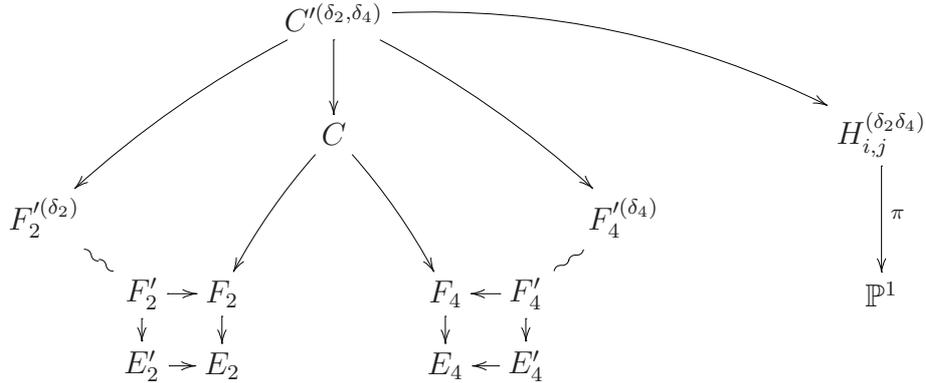
\xymatrix@R=0.7pc@C=0.8pc{
&             &            &      \ar@/_2mm/[dddlll] C'^{(\delta_2,\delta_4)} \ar@{->}[dd]    \ar@/^2mm/[dddrrr]   \ar@/^6mm/[ddrrrrrrr] &           &                           &      &       &     &   & \\
&                             &            &       &           &                           &      &    &        &   \\
&                     &            &   \ar@/_1mm/[ddl]     C      \ar@/^1mm/[ddr]                               &           &                           &   &  & & & H_{i,j}^{(\delta_2\delta_4)} \ar@{->}[dd]^\pi  \\
F_2'^{(\delta_2)}\ar@{~}[dr] &                                           &            &       &           &                          &   \ar@{~}[dl]  F_4'^{(\delta_4)}   &       &        & & \\
& F_2'  \ar@{->}[d] \ar@{->}[r]&  F_2 \ar@{->}[d]   &                                              &  F_4  \ar@{->}[d] &\ar@{->}[l] F_4' \ar@{->}[d] &  &  &  &&  \mathbb{P}^1   & \\
& E_2'  \ar@{->}[r]&  E_2    &                                              &  E_4   &\ar@{->}[l] E_4'  &      &   & & & &
}
$$

\section{An Elliptic Chabauty Argument}

Our aim in this section is to compute, for any choice of a sign
$\pm$, all the $\bQ(\sqrt{6})$-rational points in some of the
curves $H_{i,j}^{\pm } \ : \ \pm w^2=q_i(t)p_j(t)$ such that $t\in
\bQ$. We will be able to do these once we have that the jacobians
of the curves have rank 0 or 1 over $\bQ(\sqrt{6})$, a condition
coming from the Chabauty technique.

Since we only need to show this result for just one $(i,j)$, we
will do it for $(1,1)$ for both signs. This choice is not totally
arbitrary, since one can show that all the cases with $j=2$ do
have rank $2$, hence they do not fulfill the necessary conditions.

We will denote by $H^{\pm}=H^{\pm}_{1,1}$, which are the genus one
curves defined over $\bQ(\sqrt{6})$ by the equations
$$ \pm
w^2=q_1(t)p_1(t)=(5+2\sqrt{6})\,(6\,t^2+2\sqrt{6}\,t+5-6\sqrt{6})\,(6\,t^2-4\,t-1).$$

\begin{lem}\label{jacobiansH} Consider the points $P^+_{\pm}:=(1,\pm 1)\in
H^{+}(\bQ(\sqrt{6}))$ and $P^-_{\pm}:=(1/2, \pm (-2\sqrt{6} - 3)/2
)\in H^{-}(\bQ(\sqrt{6}))$. Then, the curves $H^{\pm}$ are
isomorphic over $\bQ(\sqrt{6})$ to their corresponding jacobians
$J^{\pm}:=\Jac(H^{\pm})$, which are given by the equations
$$\begin{array}{l}
J^{+} \ : \  y^2 + (-2\sqrt{6} - 10)xy + (38\sqrt{6} + 22)y = \\
 \qquad\qquad\qquad\qquad\qquad\qquad\quad x^3 + (-24\sqrt{6} - 34)x^2 + (448\sqrt{6} + 1253)x
\end{array} $$
and
$$
\begin{array}{l}
J^{-} \ : \  y^2 + (42\sqrt{6} + 422)xy + (-291822\sqrt{6} -113902)y = \\
  \qquad\qquad\qquad x^3 + (-8076\sqrt{6} - 33466)x^2 + (67635708\sqrt{6} +
141575953)x , \end{array} $$ by isomorphisms $\mu^{\pm}: H^{\pm}\to J^{\pm}$
sending the points $P^{\pm}_{+}$ to the zero point in $J^{\pm}$. Moreover,
$\mu^{\pm}(P^{\pm}_{-})=-(0,0)$.

 Finally, we have
that a point $(t,w)\in H^{\pm}(\bQ(\sqrt{6}))$ verifies that $t\in
\bQ$ if and only if $\pi^{\pm}(x,y)\in \bP^1(\bQ)$, where
$(x,y)=\mu^{\pm}(t,w)\in J^{\pm}(\bQ(\sqrt{6}))$ denotes the corresponding
image of $(t,w)$ and
$$\pi^+(x,y)=\frac {y}{2x + y} \quad \mbox{ and }\quad \pi^-(x,y)=\frac {y}{300x +
2y}.$$
\end{lem}

\begin{pf} The first part is an standard computation (see for example \cite[\S 7.2.3]{Coh}
). The inverse of the maps $\mu^{\pm}$ are given by the
maps $\nu^{\pm}$ defined as
$$
\nu^+(x,y):=\left( \frac y{2x + y},
\frac{2x^3 + (-24\sqrt{6}-34)x^2+(2\sqrt{6}+ 10)xy-y^2}{2x+ y}
\right)
$$
and
$$
\begin{array}{l}
\nu^-(x,y):=\left( \frac{y}{300x + 150y},\frac{
-2(2\sqrt{6}+3)x^3+(91160\sqrt{6}+197310)x^2+(970\sqrt{6}+1770)xy+(2\sqrt{6}+3)y^2
}{150^3(150x+ y)} \right)\!.
\end{array}
$$
The maps $\pi^{\pm}$ are the composition of the map
$\nu^{\pm}$ and the map $H^{\pm}\to \bP^1$ sending $(t,w)$ to $t$.
\end{pf}

\begin{rem}
Note that $j(J^{\pm})=140608/245\in\bQ$. Moreover,  $J^{\pm}$ is
isomorphic to the  $\pm(\sqrt{6}-3)$-twist of the elliptic curve
defined over $\bQ$ given by the Weierstrass model
$y^2=x^3+312x-3008$ (which is \verb+80640CU2+ in Cremona's
tables).
\end{rem}

\begin{lem}\label{rankJ}
The elliptic curves $J^{\pm}$ have both rank 1 over $\bQ(\sqrt{6})$.
We have
$$
J^{\pm}(\bQ(\sqrt{6}))\supset S^{\pm}= \langle T ^{\pm} ,(0,0)\rangle,
$$
where $$T^+:=(2\sqrt{6} - 7 , -16\sqrt{6} - 34) \mbox{ and }
T^-:=(-462\sqrt{6} - 2767 , 301500\sqrt{6} + 699000 )$$ are
$2$-torsion points and $(0,0)$ is of infinite order. The subgroups
$S^{\pm}$ are of finite index not divisible by any prime $<14$.
\end{lem}

\begin{pf} This is shown by an standard $2$-descent argument,
either using \verb+Magma+, \verb+Sage+ or \verb+PARI+ (see for
example \cite[\S 8.3]{Coh}). The non divisibility condition of the
index can be shown easily by proving that the point of infinite
order $(0,0)$ in both cases is not a $p$-multiple of another point
in $J^{\pm}(\bQ(\sqrt{6}))$ for any prime $p<14$.
\end{pf}

\

Now, we are under the conditions to apply the Chabauty technique.
We need to choose a prime $p$ of good reduction for $J^{\pm}$, and
also inert in $\bQ(\sqrt{6})$ (the technique can also be used for
split primes, but with a slightly distinct form, see for example
\cite{Br}). Denote by $J^{\pm}_p$ the reduction modulo $p$ of $J^{\pm}$,
which is an elliptic curve over $\bF_{p^2}:=\bF_p(\sqrt{6})$, and
by $\mbox{red}_p:J^{\pm}(\bQ(\sqrt{6})) \to J^{\pm}_p(\bF_{p^2})$ the
reduction map. Then the Elliptic Chabauty method will allow us to
bound, for each point $R$ in $J^{\pm}_p(\bF_{p^2})$, the number of
points $Q$ in $J^{\pm}(\bQ(\sqrt{6}))$ such that $\mbox{red}_p(Q)=R$
and such that $\pi^{\pm}(Q)\in\bP^1(\bQ)$. Denote this set of points
by
$$\Omega_{\pm,p}(R):=\{Q\in J^{\pm}(\bQ(\sqrt{6}))\ | \
\pi^{\pm}(Q)\in \bQ \mbox{ and } \mbox{red}_p(Q)=R\}.$$ Clearly,
we have
$$\{Q\in J^{\pm}(\bQ(\sqrt{6}))\ | \
\pi^{\pm}(Q)\in \bQ\} = \bigsqcup_{R\in J^{\pm}_p(\bF_{p^2})}
\Omega_{\pm,p}(R),$$ for any choice of an inert good reduction
prime $p$. So, if we compute these sets $\Omega_{\pm,p}(R)$ for
some choice of the prime $p$ and for all $R$, we have computed the
sets we are interested on.

We will choose the primes $p=11$ and $p=13$, depending of the sign
of the case considered.

\begin{prop}\label{badpoints} Denote by $p_{+}=11$ and $p_{-}=13$. Then
$\Omega_{\pm,p_{\pm}}(\widetilde{R})\ne \emptyset$ if and only if
$\widetilde{R}=O$ or $\widetilde{R}=-(0,0)$.
\end{prop}

\begin{pf} First of all, observe that
$\Omega_{\pm,p_{\pm}}(O)\ne \emptyset$ and
$\Omega_{\pm,p_{\pm}}(-(0,0))\ne \emptyset$ since they contain the
points $O$ and $-(0,0)$, respectively.

In order to show the remaining subsets are empty, we will use
arguments modulo $p_{\pm}^m$ for various powers of $p_{\pm}$. Denote
by $\cO$ the ring of integers of $\bQ(\sqrt{6})$, by $\cJ$ the
N{\'e}ron model of $J$ over $\cO$ and by
$\pi^{\pm}_{p_{\pm}^n}:\cJ^{\pm}_{\cO/p_{\pm}^n\cO} \to \bP^1$ the
reduction modulo $p_{\pm}^n$ of the map $\pi^{\pm}$, which is a
well-defined map of schemes over $\cO/p_{\pm}^n\cO$. Observe that
for good reduction primes $p$ as we have, the scheme
$\cJ^{\pm}_{\cO/p^n\cO}$ is an abelian scheme.

First we work modulo $p_{\pm}$. We get that $(0,0)$ has order $8$
in $J^+_{11}(\bF_{11^2})$ and has order 12 in
$J^-_{13}(\bF_{13^2})$. Since the point $(0,0)$ is not divisible
by $2$ and $3$ in $J^{\pm}(\bQ(\sqrt{6}))$ as shown in lemma
\ref{rankJ}, we get in both cases that
$\mbox{red}_{p_{\pm}}(S^{\pm})=\mbox{red}_{p_{\pm}}(J^{\pm}(\bQ(\sqrt{6})))$,
so we can work with the subgroup $S^{\pm}$.

One easily computes that the only points $R$ in
$\mbox{red}_{p_{\pm}}(S^{\pm})$ such that $\pi^{\pm}_{p_{\pm}}(R)\in
\cO/{p_{\pm}}\cO=\bF_{p_{\pm}^2}$ are, in each case, the points
$$O, -(0,0), T^++2(0,0), T^+-3(0,0) \in \mbox{red}_{p_+}(S^{+}) $$
and
$$O, -(0,0), 4(0,0),-5(0,0)\in \mbox{red}_{p_-}(S^{-}).$$

Since $\Omega_{\pm,p_{\pm}}(R)$ is obviously empty if
$\pi^{\pm}_{p_{\pm}}(R)$ is not in $\bF_{p_{\pm}^2}$, we only need to
show that $\Omega_{+,p_{+}}(R) = \emptyset$ if $R=T^++2(0,0)$ or
$T^+-3(0,0)$, and that $\Omega_{-,p_{-}}(R) = \emptyset$ if
$R=4(0,0)$ or $-5(0,0)$.

We start with the $+$ case. In this case one computes all the
points in $\cJ^+_{\cO/11^2\cO}$ which are equal to $R=T^++2(0,0)$
or to $T^+-3(0,0)$ modulo $11$ (there are $22$ of them), and then
we compute their images by $\pi^+_{11^2}$. But for any of these
$22$ points, the image is not in $\bZ/11^2\bZ \hookrightarrow
\cO/11^2\cO$, hence $\Omega_{+,p_{+}}(R) = \emptyset$ for both
points.

Finally, the $-$ case is done similarly, but one needs to work
modulo $13^3$, since all the lifts to $13^2$ do have image in
$\bZ/13^2\bZ$ with respect to the map $\pi^+_{13^2}$. When working
modulo $13^3$, the total number of points to be considered is
$2\cdot 13^2$.
\end{pf}

\begin{prop}\label{goodpoints} The sets $\Omega_{\pm,p_{\pm}}(O)$ and
$\Omega_{\pm,p_{\pm}}(-(0,0))$ contain, for any sign, only one point.
\end{prop}

\begin{pf} We use the Chabauty argument. First of all,
recall that the order of $(0,0)$ modulo $p_{\pm}$ is either
$m_+:=8$ in the $+$ case and $m_-:=12$ in the $-$ case. Hence, any
point in $\Omega_{\pm,p_{\pm}}(R)$ must be of the form
$R+n(m_{\pm}(0,0))$ for some $n\in \bZ$. Since in both cases,
$R=O$ and $R=-(0,0)$, we do have that $R$ is in
$\Omega_{\pm,p_{\pm}}(R)$, we want to show that the only solution
is $n=0$ in all cases. We will work with an argument modulo
$p_{\pm}^2$.

We will use the expression of the $z$-coordinate of the point $(x,y)$
(with respect to the given equation of $J^{\pm}$) to mean
$z:=-x/y$ (as a point in $\bP^1$).

Let us denote by $z_{\pm,p}$ the $z$-coordinate of the point
$m_{\pm}(0,0)$ modulo $p_{\pm}^2$. We get
$z_{+,11}=11-55\sqrt{6}\in \cO/11^2\cO$ and $z_{-,13}=26-39\sqrt{6} \in \cO/13^2\cO$. Because we are working
modulo $p_{\pm}^2$, we have 
$$
z\mbox{-coord}(n(m_{\pm}(0,0)))=nz_{\pm,p_{\pm}} \pmod{p_{\pm}^2\cO}
$$
(a fact that can be proved using the formal logarithm and
exponential of the elliptic curves $J^{\pm}$).

Now, we can express the function $\pi_{\pm}(P)$ at any point $P$
as a power series in the $z$-coordinate of P. We get, for the $+$
case, that
$$\pi_{+}(z)= 1 + 2z + 4z^2 + 8z^3 + O(z^4) $$
and for the $-$ case, that
$$\pi_{-}(z)=1/2 + 75z + 11250z^2 + 1687500z^3 + O(z^4).$$

First we treat the point $O$. We have that
$\pi_{\pm}(n(m_{\pm}(0,0)))$ can be expressed as a power series
$\Theta(n)$ in $n$ with coefficients in $\bQ(\sqrt{6})$. We
express this power series as
$\Theta(n)=\Theta_0(n)+\sqrt{6}\Theta_1(n)$, with $\Theta_i(n)$
now being a power series with coefficients in $\bQ$. Then
$\pi_{\pm}(n(m_{\pm}(0,0)))\in \bQ$ for some $n\in \bZ$ if and only
if $\Theta_1(n)=0$ for that $n$. Observe also that, since
$\pi_{\pm}(O)\in \bQ$, we will get that $\Theta_1(0)=0$, so
$\Theta_1(n)=j_1n+j_2n^2+j_3n^3+\cdots$. To conclude, we will use
Strassmann Theorem: if the $p_{\pm}$-adic valuation of $j_1$ is
strictly smaller that the $p_{\pm}$-adic valuation of $j_i$ for
any $i>1$, then this power series has only one zero at the
$p$-adic ring $\bZ_{p_{\pm}}$, and this zero is $n=0$. In fact,
one can easily shown that this power series verifies that the
$p_{\pm}$-adic valuation of $j_i$ is always greater or equal to
$i$, so, if we show that $j_1\not{\!\!\equiv} 0 \pmodd{p_{\pm}^2}$
we are done.

Since the $z$-coordinate of $m_{\pm}(0,0)$ is congruent to $0$
modulo $p_{\pm}$, to compute
$\pi_{\pm}(z\mbox{-coord}(n(m_{\pm}(0,0))))$ modulo $p_{\pm}^2$ we
only need the power series up to degree $1$. We get 
$$\pi_{+}(z\mbox{-coord}(n(m_{+}(0,0))))=(22-44\sqrt{6})n + 1 \pmod{11^2}$$
and
$$\pi_{-}(z\mbox{-coord}(n(m_{-}(0,0))))=
(-78-52\sqrt{6})n -84 \pmod{13^2}. $$ Hence, for the $+$ case,
$\Theta_1(n)=-44n$ modulo $11^2$, thus the valuation of $j_1$ is
$1$, and we are done. Similarly, for the $-$ case,
$\Theta_1(n)=-52n$ modulo $13^2$, and again we are done.

In order to consider the other point $-(0,0)$ one can either
compute directly $\pi_{\pm}(z\mbox{-coord}(-(0,0)+n(m_{\pm}(0,0))))$
as a power series in $n$ and apply the same type of arguments, or
observe that there is an involution in $J^{\pm}$ (as a genus 1
curve) that interchanges the points $O$ and $-(0,0)$ and preserves
the function $\pi_{\pm}$; it corresponds to the hyperelliptic
involution on $H^{\pm}$ sending $(t,w)$ to $(t,-w)$.
\end{pf}

\

Hence, by using the results just proved, we obtained finally the
following.

\begin{cor}\label{pointsH}
The only points $(t_{\pm},w)\in H^{\pm}(\bQ(\sqrt{6}))$ with $t_\pm \in
\bQ$ are the points with $t_+=1$ and with $t_-=1/2$.
\end{cor}

\section{Proof of the theorem \ref{Mthm}}

By using the results proved in the last two sections, we have
finally computed all the rational points in the curve $C$.

\begin{thm}\label{pointsC}
$C(\bQ)=\{[\pm1:\pm1:\pm1:\pm1:\pm1],[\pm1:\pm3:\pm5:\pm7:\pm9]\}$.
\end{thm}

\begin{pf}
We review here briefly the steps we followed. By using Corollary
\ref{deltaCinv}, we have that, for any point $P\in C(\bQ)$, there
exists a sign $\pm$ and an involution $\tau \in \Upsilon$ of $C$
such that $$\tau(P) \in \chi^{((\pm1,1))}(\{(t,y_1,y_2,z_1,z_2)\in
C'^{((\pm1,1))}(\bQ(\sqrt{6}))\ : \ t\in \bQ\}).$$ Moreover, the
points $[\pm1:\pm1:\pm1:\pm1:\pm1]$ and
$[\pm1:\pm3:\pm5:\pm7:\pm9]$ come, respectively, from the points
in $C'^{((1,1))}$ with $t=1$ and the points in $C'^{((-1,1))}$
with $t=1/2$.

After that, we consider the genus one quotients $H^{\pm}$ of the
curves $C'^{((\pm1,1))}$ with quotient maps defined by
$(t,y_1,y_2,z_1,z_2)\mapsto (t,y_1z_1)$. We get that the points
$(t,y_1,y_2,z_1,z_2)\in C'^{((\pm1,1))}(\bQ(\sqrt{6}))$ go to
points $(t,w)\in H^{\pm}(\bQ(\sqrt{6}))$. Finally, by Corollary
\ref{pointsH}, we have that the only points in $H^{+}$ with $t\in
\bQ$ are the ones with $t=1$, and the only points in $H^{-}$ with
$t\in \bQ$ are the ones with $t=1/2$. This proves the result.
\end{pf}

\

\begin{pfm}
If $N$ is odd, then translating by $-r-(N-1)/2$ we can suppose
$r=(N-1)/2$. Thus, $b=0$ and $f(x)=ax^2+c$ is symmetric with
respect the axis $x=0$. In the even case, we can apply the
translation by $-r-N/2$, and we can suppose that $r=-N/2$. Thus,
$b=a$ and $f(x)=a(x^2+x)+c$ is symmetric with respect the axis
$x=-1/2$. Now, the existence of $f(x)$ is equivalent to the
existence of $x_k\in \bZ$ such that $f(k)=x_k^2$,
$k=0,1,\dots,s_N$, where $s_N=(N-1)/2$ or $s_N=N/2-1$ depending on
$N$ is odd or even, respectively. Note that for $N\le 4$ is
trivial to prove that there are infinitely many non-square
quadratic satisfying the hypothesis. If $N=5$ (resp. $N=6$), then
they satisfy $3x_0^2+x_2^2=4x_1^2$ (resp. $2x_0^2+x_2^2=3x_1^2$),
that it is a conic in $\bP^2$ with infinitely rational points. If
$N=7$ (resp. $N=8$), then they satisfy $3x_0^2+x_2^2=4x_1^2$ and
$8x_0^2+x_3^2=9x_1^2$ (resp. $2x_0^2+x_2^2=3x_1^2$ and
$5x_0^2+x_3^2=6x_1^2$), which is isomorphic to the elliptic curve
$y^2=x(x-5)(x+27)$ (resp. $y^2=x(x-12)(x-15)$) and it is denoted
by \verb+30A2+ (resp. \verb+360E2+) at Cremona's table with
Mordell-Weil group isomorphic to $\bZ/2\bZ\oplus\bZ/6\bZ$ (resp.
$\bZ/2\bZ\oplus\bZ/2\bZ\oplus \bZ$). Therefore if $N\ge 7$ odd the
proof of the theorem is finished. If $N=8$, since the rank of the
underlying elliptic curve is non zero then there are infinitely
many non-square quadratic polynomials satisfying the hypothesis. The remaining
case to finish the proof is when $N\ge 10$ even. The
characterization given at section \ref{sec2}, give us that if
$N=10$ then any quadratic polynomial $f(x)\in\bZ$ satisfying the
hypothesis of the theorem correspond to a point $P\in C(\bQ)$. At
Theorem \ref{pointsC} we have proved that the unique points at
$C(\bQ)$ are $[\pm1:\pm1:\pm1:\pm1:\pm1]$,
$[\pm1:\pm3:\pm5:\pm7:\pm9]\}$ that correspond to the constant
polynomials and to $f(x)=(2x+1)^2$ respectively.
\end{pfm}

\

{\bf Data:} All the \verb+Magma+ sources are available from the first author webpage.

\

\end{document}